\documentclass{birkjour}
\usepackage{amscd}
\usepackage{amssymb}
 \newtheorem{thm}{Theorem}[section]
 \newtheorem{cor}[thm]{Corollary}
 
 \newtheorem{prop}[thm]{Proposition}
 \theoremstyle{definition}
 
 \theoremstyle{remark}

 \numberwithin{equation}{section}
\usepackage[final]{hyperref}
\hypersetup{unicode= false, colorlinks=true, linkcolor=blue,
anchorcolor=blue, citecolor=red, filecolor=red, menucolor=blue,
pagecolor=blue, urlcolor=blue}
\usepackage[none]{hyphenat}
\renewcommand{\theenumi}{\roman{enumi}}
\renewcommand{\labelenumi}{(\theenumi)}

\newcommand{\CC}{\mathbb {C}}

\title[ Volterra type  and weighted
 composition operators on  Fock spaces ] {Volterra
type and weighted composition operators on  weighted Fock spaces}
\author [Tesfa  Mengestie]{Tesfa  Mengestie }
\address{Department of Mathematical Sciences\\
Norwegian University of Science and Technology (NTNU)\\
 NO- 7491 Trondheim, Norway}
\email{tesfantnu@gmail.com}
\begin{document}
\begin{abstract}Bounded and compact  product of  Volterra type integral and
composition operators acting between weighted  Fock spaces are
described. We also estimate the norms  of these
operators in terms of Berezin type integral transforms on the
complex plan $\CC$. All our results are valid for weighted
composition operators acting on the class of  Fock spaces
considered under appropriate  interpretation of the weights.
\end{abstract}
 \maketitle
\section{Introduction}
Given a  holomorphic function $g$ on the complex plane
$\CC,$ we define the induced  Volterra type integral operator  by
 \begin{equation*}
 V_gf(z)= \int_0^z f(w)g'(w) dw.
\end{equation*}
Questions about boundedness, compactness, and other  operator
theoretic properties of $V_g$ expressed in terms of function
theoretic conditions on the symbol  $g$ have been a subject of high interest
since introduced by Pommerenke \cite{Pom} in 1997. The operator has especially received considerable
 attentions following the  works of Aleman and Siskakis \cite{Alsi1,Alsi2} on Hardy and
Bergmann spaces.  For more information on the subject,  we refer to the surveys in
\cite{Alman,Si} and the related references therein.

 Each entire function  $\psi$ induces a  composition operator defined by
  $C_\psi f=f(\psi)$ on the space of holomorphic functions on $\CC.$  Using $V_g$ and $C_\psi$, we define the
product of Volterra type integral and composition operators induced
by the pair of symbols $(g,\psi)$ as
\begin{equation*}
\label{define} V_{(g, \psi)}f(z)= (V_{g} \circ C_\psi )f(z)=  \int_0^z f(\psi(w))g'(w)dw
\end{equation*} and
\begin{equation*}
\label{define} C_{(\psi, g)}f(z)= (C_\psi \circ V_g )f(z)=  \int_0^{\psi(z)} f(w))g'(w)dw.
\end{equation*}
If $\psi(z)=z,$ then these operators are just  the usual Volterra
type integral operators $V_g.$ As will be indicated latter, studying
the  operators $V_{(g, \psi)} $ reduces to studying  the composition operator
$C_\psi$ between Fock spaces when $|g'(z)|/(1+|z|)$ behaves like a constant function. Recent years have seen
 a lot of work on  operators of these
kinds acting on different space of holomorphic functions;  for example in \cite{LIS, SLI, TM, AS, XZHU,EW}.  Inspired by the works in \cite{ZZHH, ZZH},
 S. Ueki \cite{UEKI}  characterized the bounded  weighted composition
operators on the classical Fock space $\mathcal{F}_1^2$ in terms of Berezin type integral transforms. Recently, the author
obtained  similar  results  for the operator $V_{(g,\psi)}$ when it
acts  between the  Fock spaces $\mathcal{F}_\alpha^p$ and
$\mathcal{F}_\alpha^q$ whenever both exponents $p$ and $ q$ are
finite \cite{TM}. The purpose of the present work is partly to continue that line of research
and fill the gap when  at least one of the Fock spaces in question is the growth space, i.e,  the case when one  of the exponents is at infinity. We, in addition, establish
some interesting results on  the mapping  properties of  $C_{(\psi,g)}$ and weighted composition operators on Fock spaces.
 Our work unifies and  extends a number of results on Volterra type integral  and
 weighted composition operators; for example in \cite{CMS,Olivia,TM,SS, UEKI,UEKI2}. As a main tool in proving  some of our results, we will also
characterize  the $(\infty, p)$ Fock--Carleson measures which  is an interesting result in its own right.

For $ \alpha >0, $ the classical weighted Fock space $\mathcal{F}_\alpha^p$ consists of
entire functions $f$ for which
\begin{equation*}
\|f\|_{(p,\alpha)}^p=  \frac{\alpha p}{2\pi}\int_{\CC} |f(z)|^p
e^{-\frac{\alpha p}{2}|z|^2} dm(z) <\infty
\end{equation*} for $ 0 < p <\infty$ where $dm$ is the usual Lebesgue measure on $\CC$, and for $p= \infty,$
\begin{equation*}
\|f\|_{(\infty,\alpha)}= \sup_{z\in \CC}
|f(z)|e^{-\frac{\alpha}{2}|z|^2} <\infty.
\end{equation*}  The space  $\mathcal{F}_\alpha^2$ is, in particular,  a reproducing kernel Hilbert space with
kernel and normalized  kernel  functions  respectively
$K_{(w,\alpha)}(z)= e^{\alpha \langle z, w\rangle}$ and
$k_{(w,\alpha)}(z)= e^{\alpha\langle z, w\rangle}-\alpha|w|^2/2.$

 Our  results are expressed in terms of the functions
 \begin{align*}
 B^\infty_{(\psi,\alpha)} (|g|)(z)
 = \frac{|g'(z)|}{1+|z|}e^{\frac{\alpha}{2}\big(|\psi(z)|^2-|z|^2\big)},\  \
 M^\infty_{(\psi,\alpha)} (|g|)(z)
 = \frac{|g'(\psi(z))|e^{\frac{\alpha}{2}\big(|\psi(z)|^2-|z|^2\big)}}{(1+|z|)(|\psi'(z)|^{-1})}
 \end{align*}
 and Berezin type integral transforms
 \begin{equation*}
B_{(\psi,\alpha)}(|g|^p)(w)= \int_{\CC} \frac{|k_{(w,\alpha)}(\psi(z))g'(z)|^p}{\big((1+|z|)e^{\frac{\alpha}{2}|z|^2}\big)^p} dm(z) \ \
\end{equation*} and
\begin{equation*}
M_{(\psi,\alpha)}(|g|^p)(w)= \int_{\CC} \frac{|k_{(w,\alpha)}(\psi(z))g'(\psi(z)) \psi'(z)|^p}{\big((1+|z|)e^{\frac{\alpha}{2}|z|^2}\big)^p} dm(z).
\end{equation*}
A word on notation: The notation $U(z)\lesssim V(z)$ (or equivalently $V(z)\gtrsim
U(z)$) means that there is a constant $C$ such that $U(z)\leq CV(z)$
holds for all $z$ in the set in question, which may be a Hilbert
space or  a set of complex numbers. We write $U(z)\simeq V(z)$ if
both $U(z)\lesssim V(z)$ and $V(z)\lesssim U(z)$.
\section{The main results}
\subsection{Bounded and compact $V_{(g, \psi)}$ and $C_{(\psi, g)}$ .}
In this subsection, we  formulate the  main results on the bounded and compact mapping properties of
the operators $V_{(g, \psi)}$ and $C_{(\psi, g)}$.
 \begin{thm}\label{bounded1}
 Let $0<p\leq \infty$ and $\psi$ and $g$ be entire functions on  $\CC$.   Then
 \begin{enumerate}
  \item   $V_{(g, \psi)}: \mathcal{F}_\alpha^p \to \mathcal{F}_\alpha^\infty$ is
 bounded if and
 only if  $B^\infty_{(\psi,\alpha)} (|g|)\in L^\infty( dm).$  Moreover, if $V_{(g, \psi)}$ is
 bounded, then its norm is estimated by
\begin{equation}
\label{normbound} \|V_{(g, \psi)}\|\simeq \sup_{z\in \CC}
B^\infty_{(\psi,\alpha)} (|g|)(z).\end{equation}
\item  A bounded map  $V_{(g,
\psi)}: \mathcal{F}_\alpha^p \to \mathcal{F}_\alpha^\infty$ is compact if and only if
 \begin{equation}
 \label{comp}
 \lim_{|\psi(z)|\to \infty} B^\infty_{(\psi,\alpha)} (|g|)(z)=
 0.
 \end{equation}
 \item   $C_{(\psi, g)}: \mathcal{F}_\alpha^p \to \mathcal{F}_\alpha^\infty$ is
 bounded if and
 only if  $M^\infty_{(\psi,\alpha)} (|g|)\in L^\infty(dm).$  Moreover, if $C_{(\psi, g)}$ is
 bounded, then its norm is estimated by
\begin{equation*}
\label{normbound1} \|C_{(\psi,g)}\|\simeq \sup_{z\in \CC}
M^\infty_{(\psi,\alpha)} (|g|)(z).\end{equation*}
\item  A bounded map  $C_{(\psi,g)}: \mathcal{F}_\alpha^p \to \mathcal{F}_\alpha^\infty$ is compact if and only if
 \begin{equation*}
 \label{comp1}
 \lim_{|\psi(z)|\to \infty} M^\infty_{(\psi,\alpha)} (|g|)(z)=
 0.
 \end{equation*}
 \end{enumerate}
 \end{thm}
It is interesting to note that both the conditions in (i) and (ii)
are independent of the
 exponent $p.$ It follows that  if there exists a $p>0$ for which  $ V_{(g,
\psi)}$ is bounded (compact) from $ \mathcal{F}_\alpha^p $ to $
\mathcal{F}_\alpha^\infty$, then it is also  bounded (compact) for every other
$p$. This phenomena holds true for the operator $C_{(\psi,g)}$ as well.

As remarked earlier, when $\psi(z)=z,$ the operators $ V_{(g, \psi)}$ and $C_{(\psi,g)}$ reduce to the Volterra type
integral operator $V_g$. Theorem~\ref{bounded1} ensures that the symbol $g$   can only  grow as a power function of
degree not exceeding  respectively $2$  and  $1$  to induce a bounded
and compact  $V_g$. We formulate this  as follows.
 \begin{cor}\label{secondcor} Let $0<p\leq \infty$.
   Then  $V_{g}: \mathcal{F}_\alpha^p \to \mathcal{F}_\alpha^\infty$
 is
 \begin{enumerate}
 \item bounded if and only if $g(z)= az^2+bz+c ,$ for some $a, b,c \in
 \CC$.
 \item  compact if and only if $g(z)= az+b.$
 \end{enumerate}
 \end{cor} The corollary extends similar representations obtained in \cite{Olivia,TM}
 for the symbol $g$. \\
 Mapping $\mathcal{F}_\alpha^\infty$ into
 smaller spaces $\mathcal{F}_\alpha^p$ gives the next stronger condition.
 \begin{thm}\label{bounded3}
Let $0<p< \infty, $ and $\psi$ and $g$ be entire functions on  $\CC$.  Then  the
following are equivalent.
\begin{enumerate}
\item $V_{(g,\psi)}:\mathcal{F}_\alpha^\infty \to \mathcal{F}_\alpha^p$ is
  bounded;
  \item $V_{(g,\psi)}: \mathcal{F}_\alpha^\infty \to \mathcal{F}_\alpha^p$ is compact;
  \item  $B_{(\psi,\alpha)} (|g|^p) \in L( dm).$  Moreover,  if
  $V_{(g,\psi)}$ is bounded, then
  \begin{equation}
  \label{normless}
  \|V_{(g,\psi)}\|^{p} \simeq \int_{\CC}B_{(\psi,\alpha)}
  (|g|^p)(w)dm(w).
  \end{equation}
  \end{enumerate}
 \end{thm}
 \begin{thm}\label{bounded4}
Let $0<p< \infty$, and  $\psi$ and $g$ be entire functions on  $\CC$.  Then  the
following are equivalent.
\begin{enumerate}
\item $C_{(\psi,g)}: \mathcal{F}_\alpha^\infty \to \mathcal{F}_\alpha^p$ is
  bounded;
  \item $C_{(\psi,g)}: \mathcal{F}_\alpha^\infty \to \mathcal{F}_\alpha^p$ is compact;
  \item  $M_{(\psi,\alpha)} (|g|^p) \in L(dm).$  Moreover if
  $C_{(\psi,g)}$ is bounded, then
  \begin{equation}
  \label{normless1}
  \|C_{(\psi,g)}\|^p \simeq \int_{\CC}M_{(\psi,\alpha)}
  (|g|^p)(w)dm(w).
  \end{equation}
  \end{enumerate}
 \end{thm}
 For the operators $V_g$ and $C_\psi$, we in particular get the following
 explicit expressions for the inducing symbols.
 \begin{cor}\label{fourthcor}
 Let $0<p< \infty$ and  $\psi$ and $g$ be entire functions on  $\CC$.  Then
 \begin{enumerate}
 \item  $V_{g}: \mathcal{F}_\alpha^\infty \to \mathcal{F}_\alpha^p$ is
  bounded (compact)  if and
 only if    $g(z)= az+b$ and $p>2$ where $a$ and $b$ are in $\CC.$
 \item $C_\psi:\mathcal{F}_\alpha^\infty \to \mathcal{F}_\alpha^p$ is bounded (compact) if and only if
$\psi(z)= az+b,  |a|<1.$
 \end{enumerate}
 \end{cor}
 Part (ii) extends the results in \cite{CMS} where
 similar conditions  are  given for compact and bounded $C_\psi:\mathcal{F}_\alpha^p \to
 \mathcal{F}_\alpha^q$ when  $0<p\leq q<\infty$ and in \cite{TM} when $ 0<q< p<\infty$.
 \begin{proof}
(i) We start by proving the necessity of the statements in (i). We first observe that
 by applying the explicit expression for the normalized reproducing kernels, we have
 \begin{eqnarray}
 \label{con}
B_{(\psi,\alpha)}(|g|^p)(w)&=& \int_{\CC} \frac{|k_{(w,\alpha)}(\psi(z))g'(z)|^p}{\big((1+|z|)e^{\frac{\alpha}{2}|z|^2}\big)^p} dm(z)\nonumber\\
&=&\int_{\CC} e^{\frac{p\alpha}{2}\big(
2\Re \langle \psi(z), \
w\rangle-|z|^2-|w|^2\big)}\frac{|g'(z)|^p}{(1+|z|)^p} dm(z)\nonumber\\
&=&\int_{\CC} e^{\frac{p\alpha}{2}\big(
|\psi(z)|^2-|z|^2-|w-\psi(z)|^2\big)}\frac{|g'(z)|^p}{(1+|z|)^p} dm(z)
\end{eqnarray} where the last equality follows by completing the square on the inner product. If
$D(w,1)=\{z\in \CC: |z-w|< 1\},$ then  $1+|z| \simeq 1+|w|$ for each $z\in D(w,1) $. It follows that
\begin{align*}
\int_{\CC}B_{(\psi,\alpha)}(|g|^p)(w)dw\geq\int_{\CC} \int_{D(w,1)}
\bigg|\frac{k_{(w,\alpha)}(z)g'(z)}{(1+|z|)}\bigg|^p e^{-\frac{p\alpha}{2}|z|^2}dm(z)dm(w)\\
\geq \int_{\CC} \int_{D(w,1)}\bigg|\frac{g'(w)}{1+|w|}\bigg|^{p} e^{-\frac{p\alpha}{2}|z-w|^2} dm(z)dm(w)\\
\gtrsim \int_{\CC}\bigg| \frac{g'(w)}{1+|w|}\bigg|^{p} dm(w).
\end{align*}The
 desired restrictions on $g$ and $p$
follow once we apply Theorem~\ref{bounded3}. On the other hand, the sufficiency of the condition is
 immediate because by Theorem~\ref{bounded3} again and Fubini's  Theorem
\begin{align*}
\int_{\CC}B_{(\psi,\alpha)}(|g|^p)(w)dm(w) \simeq  \int_{\CC}
\int_{\CC} \frac{e^{\frac{-p\alpha}{2}|w-z|^2}}{(1+|z|)^p} dm(z) dm(w)
\simeq \int_{\CC} \frac{dm(z)}{(1+|z|)^{p}} <\infty
\end{align*} for each $p>2.$\\
(ii) Applying once  again Theorem~\ref{bounded3} and Fubini's  Theorem, we obtain
\begin{align*}
 \int_{\CC}\int_{\CC} \bigg| \frac{k_{(w,\alpha)}(\psi(z))}{
e^{\frac{\alpha}{2}|z|^2}}\bigg|^p
dm(z)dw(w)= \int_{\CC} \int_{\CC}e^{\frac{p
\alpha}{2}\big(|\psi(z)|^2 -|z|^2-|\psi(z)-w|^2\big)}
dm(w) dm(z)\\
  \simeq\int_{\CC}
e^{\frac{p\alpha}{2}(|\psi(z)|^2 -|z|^2)} dm(z) <\infty \ \ \ \ \ \ \ \ \ \ \ \ \
\end{align*} if and only if $\psi$ has a linear form as described
above  with $|a|<1.$
\end{proof}
\subsection{ The $(\infty, p)$ Fock--Carleson measures.}  For
$0<p<\infty$, we call  a nonnegative measure $\mu$ on $\CC$ an
$(\infty, p) $ Fock--Carleson measure if   the  inequality
\begin{equation*}
\bigg(\int_{\CC} |f(z)|^p e^{\frac{-p\alpha}{2}|z|^2}
d\mu(z)\bigg)^{1/p}\lesssim \| f\|_{(\infty,\alpha)}
\end{equation*} holds for all $f$ in $\mathcal{F}_\alpha^\infty$.
 In other words,  $\mu$ is an $(\infty, p) $ Fock--Carleson measure if and only if
the embedding map $I_\mu: \mathcal{F}_{\alpha}^\infty\to L^p(d\mu)$ is bounded where $d\mu(z)= e^{- \frac{\alpha p}{2}|z|^2}dV(z)$.  We
 call $\mu$   an $(\infty, p) $  vanishing Fock--Carleson measure when such an embedding map is compact. When both exponents are
finite,  the $(p,q)$ Fock--Carleson  measures were described in \cite{ZHXL} in terms of
the parameters t-Berezin transform of $\mu,$
\begin{equation*}
\widetilde{\mu_{(t,\alpha)}}(w)= \int_{\CC} |k_{(w,\alpha)}(z)|^{t}
e^{-\frac{\alpha t}{2}|z|^2}d\mu(z),
\end{equation*} and  $D(z, \delta)$ , the  Euclidean disc
centered at $z$ and radius $\delta>0.$  We will use  these  parameters to describe the
$(\infty,p)$ Fock--Carleson measures as well.

We first remind a classical covering  result  that there exists a
 sequence of points $z_k$ such that for any $r>0,$ the discs $ D(z_k,
r/2)$ cover $\CC$ and  the discs $D(z_k, r/4)$ are pairwise
disjoint. From these conditions it follows that for any $\delta
>0,$ there  exists  a positive integer $N_{\max}= N_{\max}(\delta, r)$ such that
every point in  $\CC$ belongs to at most $N_{\max}$ of the discs. Then,  we
prove the following result which  is interesting by its own right  apart from being needed to prove part of
Theorem~\ref{bounded3} and Theorem~\ref{bounded4} in the next subsection.
\begin{prop}\label{Carleson}
Let $0<p<\infty, $ and $ \mu >0 $ be a measure on $\CC.$ Then the
following statements are equivalent.
\begin{enumerate}
\item $\mu$ is an $(\infty, p)$ Fock--Carleson measure;
\item $\mu$ is an $(\infty, p)$ vanishing Fock--Carleson measure;
\item $\widetilde{\mu_{(t,\alpha)}} \in L(dm)$ for some (or any ) $t>0$;
 \item $\mu$ is a finite measure on $\CC$;
\item $\mu(D(., r)) \in L(dm)$ for some (or any)
$r>0$;
\item $ \mu\big(D(z_k,r))\in \ell^1$ for some (or any)
$r>0$.  Moreover, \begin{equation} \label{fourequal} \|I_\mu\|^p\simeq
\|\widetilde{\mu_{(t,\alpha)}}\|_{L( dm)}\simeq \|\mu(D(.,
\delta))\|_{L( dm)}\simeq
\|\mu\big(D(z_k,r))\|_{\ell^1}.\end{equation}
\end{enumerate}
\end{prop}
 \emph{Proof.} The
equivalences of the statements in  (iii), (v) and (vi) can be found in  Lemma~ 2.3 of \cite{ZHXL}. The proof of
(i) implies (ii)  follows from a simple modification of part of the arguments
used in the proof of Theorem~3.3 in there again. Thus, we shall
prove (i) implies (vi), (vi)
implies (i), and (iii) implies (iv). For the first implication, we need  to find a test function $f_0$ in $F_\alpha^\infty$
that would lead us to the desired conclusion. For this, we follow the classical  Luecking's
approach in \cite{DL}. First,  we observe that for each  $c_j$ in
$\ell^\infty,$  by Theorem~8.2 in \cite{SJR}, we have that
 \begin{equation*}
 f=\sum_{j=1}^\infty c_j k_{(z_j, \alpha)} \in F_\alpha^\infty \ \text{and} \
 \|f\|_{(\infty,\alpha)} \lesssim \|(c_j)\|_{\ell^\infty}. \end{equation*}
Taking  $c_k=1$ for all $k$  does the job  for our case here. Thus we set our test function
as  $$f_0=\sum_{j=1}^\infty  k_{(z_j, \alpha)}.$$ Since $\mu$ is an
$(\infty, p)$ Fock--Carleson measure,
\begin{equation*}
\int_{\CC} |f_0(z)|^p e^{-\frac{\alpha p}{2}|z|^2} d\mu(z) \leq
\|I_\mu\|^p \|f_0\|_{(\infty,\alpha)}^p \lesssim \|I_\mu\|^p.
\end{equation*}
If $r_j's(t)$ are the Rademacher sequence of functions on $[0,1]$
chosen  as in \cite{DL}, then Khinchine's inequality yields
\begin{equation}
\Bigg(\sum_{j=1}^\infty
|k_{(z_j, \alpha)}(z)|^2\Bigg)^{p/2} \lesssim \int_0^1
\bigg|\sum_{j=1}^\infty r_j(t)k_{(z_j, \alpha)}(z)\bigg|^pdt.
\label{Khinchine}
\end{equation}
Note that if the $r_j(t)$ are chosen as refereed above, then
$r_j(t)$ belongs to $ \ell^\infty$  and  $\|(r_j(t))\|_{\ell^\infty}
\simeq 1$. Moreover, it holds that
\begin{equation*}
\sum_{j=1}^\infty r_j(t)k_{(z_j, \alpha)} \in
F_\alpha^\infty, \ \text{with}\
 \Big\|\sum_{j=1}^\infty r_j(t)k_{(z_j, \alpha)} \Big\|_{(\infty,\alpha)} \simeq
1.\end{equation*}
Making use of first \eqref{Khinchine} and subsequently Fubini's
theorem, we obtain
\begin{align}
\int_{\CC} \Bigg( \sum_{k=1}^\infty
|k_{(z_j, \alpha)}|^2\Bigg)^{\frac{p}{2}} d\mu(z) \lesssim \int_{\CC}
\bigg(\int_0^1 \bigg| \sum_{j=1}^\infty
r_j(t)k_{(z_j,\alpha)}(z)\bigg|^p
dt\bigg) d\mu(z)\nonumber\\
=\int_0^1\Bigg(\int_{\CC}\bigg| \sum_{j=1}^\infty
r_j(t)k_{(z_j,\alpha)}(z)\bigg|^p d\mu(z)\Bigg)dt\lesssim\|I_\mu\|^p .
 \label{combine}
\end{align}
On the other hand,
\begin{equation*}
\sum_{j=1}^\infty \mu(D(z_j,r))=
\int_{\CC}\sum_{j=1}^\infty  \chi_{D(z_j,r)}(z)
d\mu(z)\leq \int_{\CC}\Bigg(\sum_{j=1}^\infty
\chi_{D(z_j,r)}(z)\Bigg)^{p/2}d\mu(z) \label{pbig}
\end{equation*}which follows by H\"{o}lder's inequality for $p<2$ and general fact for $p\geq 2$ as all the terms are positive.
The last integral above is bounded by
\begin{equation*}
 \int_{\CC}\Bigg(\sum_{j=1}^\infty e^{\alpha r^2}|k_{(z_j, \alpha)}|^2\Bigg)^{p/2}d\mu(z)\lesssim
\int_{\CC}\Bigg(\sum_{j=1}^\infty |k_{(z_j, \alpha)}|^2\Bigg)^{p/2}d\mu(z) .
\end{equation*}This combined with \eqref{combine} gives
\begin{equation}
\label{oneside} \sum_{j=1}^\infty \mu(D(z_j,r)) \lesssim
\|I_\mu\|^p.
\end{equation}
To prove that (i) follows from  (vi),  by Lemma~1 of \cite{JIKZ}, we  observe that
 \begin{equation*}
\label{book}
\sup_{z\in D(z_k,r)} |f(z)e^{-\frac{\alpha}{2}|z|^2 }|^p
\lesssim \int_{D(z_k,2r)}|f(z)e^{-\frac{\alpha}{2}|z|^2}|^p dm(z)
\end{equation*} for each  $k$. Thus, we get
\begin{align}
\int_{\CC}|f(z)e^{-\frac{\alpha}{2}|z|^2}|^p d\mu(z)
\leq \sum_{k=1}^\infty\int_{D(z_k,r)}|f(z)e^{-\frac{\alpha}{2}|z|^2}|^p d\mu(z)\ \ \ \ \ \ \ \ \ \ \ \ \ \ \ \ \ \ \nonumber\\
\lesssim  \|\mu(D(z_k,r))\|_{\ell^1}\sum_{k=1}^\infty\int_{D(z_k,2r)}|f(z)e^{-\frac{\alpha}{2}|z|^2}|^p dm(z)\nonumber\\
 \lesssim \|\mu(D(z_k,r))\|_{\ell^1}\|f\|_{(\infty,\alpha)}^p\nonumber
\label{book1}
\end{align} from which the estimate
\begin{equation}
\label{normpart} \|I_u\|^p \lesssim \|
\mu(D(z_k,r))\|_{\ell^1}
\end{equation} holds. Now we combine \eqref{normpart}, \eqref{oneside} and the estimates in  Theorem~3.3 of
\cite{ZHXL}  to get all the remaining norm estimates in
\eqref{fourequal}.\\
We remain to show (iii) implies (iv). But this  global  geometric condition follows
 when we in particular set  $t=1$. Because by Fubini's theorem, we may have
 \begin{align*}
 \int_{\CC} \widetilde{\mu}_{(1,\alpha)}(z) dm(z)=\int_{\CC}   \int_{\CC} e^{\frac{\alpha}{2} |\langle w, z\rangle|^2 -\frac{\alpha}{2}
 |z|^2-\frac{\alpha}{2}
 |w|^2}d\mu(w)dm(z)  \ \ \  \nonumber\\
 =\int_{\CC} \Bigg(\int_{\CC} e^{-\frac{\alpha}{2}
 |z-w|^2} dm(z)\Bigg)d\mu(w) \simeq \mu (\CC).
 \label{finite}
 \end{align*}
 \subsection{Proof of the main results.}
Using a result of H. Cho and K. Zhu \cite {HRKZ}, recently
Constantin \cite{Olivia} proved that \begin{equation} \label{olivia}
\int_{\CC} |f(z)|^p e^{-p\alpha|z|^2/2}dm(z) \simeq |f(0)|^p
+\int_{\CC} \frac{|f'(z)|^p}{(1+|z|)^p}e^{-p\alpha|z|^2/2}dm(z)
\end{equation} holds for each entire function $f$ and
$0<p<\infty.$  From a
simple variant of the arguments used in the proof of  \eqref{olivia}, we also  conclude
that
\begin{equation}
\label{compare}
 \sup_{z\in \CC} |f(z)| e^{\frac{-\alpha}{2}|z|^2}  \simeq |f(0)|+
 \sup_{z\in \CC}  \frac{|f'(z)| }{1+|z| }e^{\frac{-\alpha}{2}|z|^2}
\end{equation} for each entire function $f$. Observe that
\eqref{olivia} and \eqref{compare} describe the Fock spaces in terms of
 derivatives. Such a description  plays an important role in our  further analysis. In fact, it  makes
 our analysis easier by eliminating  the integral that
 arises from the Volterra type integral operator $V_g$.
 By Lemma~1 of \cite{SS}, we also have the pointwise estimate
 \begin{equation}
 \label{point}
 |f(z)|e^{-\frac{\alpha}{2}|z|^2} \leq \|f\|_{(p,\alpha)}
 \end{equation} for each point $z $ in $\CC, \ 0<p\leq \infty$ and $f$ in $ F_\alpha^p$.  From this along with  \eqref{compare}  it follows that
  \begin{equation*}
 \|V_{(g,\psi)}f\|_{(\infty,\alpha)} \simeq \sup_{z\in \CC}  \frac{|g'(z)| }{1+|z| }|f(\psi(z))|e^{\frac{-\alpha}{2}|z|^2} \leq \|f\|_{(p,\alpha)} \sup_{z\in \CC}B^\infty_{(\psi,\alpha)} (|g|)(z)
 \end{equation*} from which the sufficiency of  the statement in part (i) of
 Theorem~\ref{bounded1} and the estimate
 \begin{equation}
 \label{one}
 \|V_{(g,\psi)}\|\lesssim \sup_{z\in \CC} B^\infty_{(\psi,\alpha)}(|g|)(z)
 \end{equation} follow. To prove its necessity, we apply $V_{(g,\psi)}$ to the normalized  kernel function $k_{(w,\alpha)}$ and estimate
  the resulting function $(\infty,\alpha)$ norm using  \eqref{compare} as
 \begin{eqnarray*}
 \|V_{(g,\psi)}k_{(w,\alpha)}\|_{(\infty,\alpha)} &\simeq& \sup_{z\in \CC}  \frac{|g'(z)| }{1+|z| }|k_{(w,\alpha)}(\psi(z))|e^{\frac{-\alpha}{2}|z|^2}\nonumber\\
 &\geq &\frac{|g'(z)| }{1+|z| }|k_{(w,\alpha)}(\psi(z))|e^{\frac{-\alpha}{2}|z|^2}
 \end{eqnarray*}  for all $z, w\in \CC.$ In particular,  when we set $w= \psi(z)$, we get
 \begin{equation*}
 \|V_{(g,\psi)}k_{(w,\alpha)}\|_{(\infty,\alpha)} \gtrsim B^\infty_{(\psi,\alpha)} (|g|)(z)
 \end{equation*} from which the reverse inequality in \eqref{one} holds.\\
 To prove the second part of the theorem, we extend a technique used  in
  \cite{SS,UEKI2}. We first assume that $ V_{(g,\psi)}$ is compact. The  sequence  $k_{(w, \alpha)}$ converges to zero as
$|w| \to \infty,$ and the convergence is uniform on compact subset of $\CC.$  We further assume that there exists sequence of points
$z_j \in \CC$ such that $|\psi(z_j)| \to \infty$ as $j\to \infty.$ If such a sequence does not exist, then necessity holds trivially. It follows from compactness of  $V_{(g,\psi)}$ that
\begin{equation}
\limsup_{j\to \infty} B^\infty_{(\psi,\alpha)}(z)  \leq \limsup_{j\to \infty}\| V_{(g,\psi)} k_{(w, \alpha)} \|_{(\infty,\alpha)}=0
\end{equation} from which \eqref{comp} follows.

We now assume that condition \eqref{comp} holds, and proceed to show that $V_{(g,\psi)} $ is a compact map. The function $f(z)= 1$ belongs to $\mathcal{F}_{\alpha}^p$ for all $p>0.$  It follows that by  boundedness,
\begin{equation}
\label{00}
\sup_{z\in \CC} \frac{|g'(z)|}{ 1+|z|} e^{-\frac{\alpha}{2}|z|^2}<\infty.\end{equation}
Let $f_j$ be a sequence of functions in $\mathcal{F}_{\alpha}^p$ such that $\sup_m
\|f_j\|_{(p,\alpha)}<\infty$ and $f_j$ converges uniformly to
zero on compact subsets of $\CC$ as $j\to \infty.$  For each $\epsilon >0,$ condition \eqref{comp} implies that there exists
a positive $N_1$ such that
$ B^\infty_{(\psi,\alpha)}(z) <\epsilon
$ for all $|\psi(z)| > N_1.$ From this together with \eqref{point}, we obtain
\begin{equation}
\label{tmo}
 \frac{|g'(z)|}{ 1+|z|}| f_j(\psi(z)) | e^{-\frac{\alpha}{2}|z|^2}
   \leq \|f_j\|_{(p,\alpha)}\frac{|g'(z)|}{ 1+|z|} e^{\frac{\alpha}{2}|\psi(z)|^2-\frac{\alpha}{2}|z|^2}\lesssim \epsilon
 \end{equation} for all $|\psi(z)| > N_1$ and all $j.$ On the other hand if $|\psi(z)|\leq N_1,$ then applying \eqref{00} it is easily seen that
\begin{equation}
\label{tm1}
  \frac{|g'(z)|}{ 1+|z|}| f_j(\psi(z)) | e^{-\frac{\alpha}{2}|z|^2}
  \lesssim \sup_{z:|\psi(z)|\leq N_1 } |f_j(\psi(z))| \to 0
   \end{equation} as $j\to \infty$. Then we apply \eqref{tmo},\eqref{tm1}  and \eqref{compare}
   to arrive at the desired conclusion.\\
   (iii) We first note that
 $(C_{(\psi,g)}f(z))'= g'(\psi(z)) \psi'(z) f(\psi(z)).$   This together with
 \eqref{compare} implies
 \begin{align*}
 \|C_{(\psi,g)}f\|_{(\infty,\alpha)} \simeq \sup_{z\in \CC}  \frac{|g'(\psi(z))\psi'(z)| }{1+|z| }|f(\psi(z))|e^{\frac{-\alpha}{2}|z|^2}\\ \lesssim
  \|f\|_{(p,\alpha)} \sup_{z\in \CC}\frac{|g'(\psi(z))\psi'(z)| }{1+|z|}e^{\frac{\alpha}{2}|\psi(z)|^2-\frac{\alpha}{2}|z|^2} \\
  =  \|f\|_{(p,\alpha)} \sup_{z\in \CC} M^\infty_{(\psi,\alpha)} (|g|)(z)
 \end{align*} where the second inequality follows by \eqref{point}. From this the sufficiency part of (iii)
 and the estimate
 \begin{equation}
 \label{onee}
 \|C_{(\psi,g)}\| \lesssim \sup_{z\in \CC} M^\infty_{(\psi,\alpha)} (|g|)(z)
 \end{equation} follow. To prove the necessity, we apply
 $C_{(\psi,g)}$ to the normalized kernel function again and observe that
 \begin{align*}
 \|C_{(\psi,g)}k_{(w,\alpha)}\|_{(\infty,\alpha)} \gtrsim \frac{|g'(\psi(z))\psi'(z)| }{1+|z| }|k_{(w,\alpha)}(\psi(z))|e^{\frac{-\alpha}{2}|z|^2}
 \end{align*}  for all $z, w\in \CC.$ In particular when we set $w= \psi(z)$, we again get
 \begin{equation*}
 \|C_{(\psi,g)}k_{(w,\alpha)}\|_{(\infty,\alpha)} \gtrsim \sup_{z\in \CC} M^\infty_{(\psi,\alpha)} (|g|)(z)
 \end{equation*} from which the reverse inequality in \eqref{onee} also follows.\\
   The proof of part (iv) of the theorem is very similar to the proof of part (ii). We thus omit it.\\
\emph{Proof of Theorem~\ref{bounded3}.}
 For each $p>0,$ we set $\mu_{(p,\alpha)}$ to be
the positive pull back  measure on $\CC$ defined by \begin{equation*}
\mu_{(p,\alpha)}(E)=\int_{\psi^{-1}(E)}
\frac{|g'(z)|^p}{(1+|z|)^p} e^{-\frac{p\alpha}{2}|z|^2}
dm(z)\end{equation*} for every Borel  subset $E$ of $\CC$.  By \eqref{olivia} and substitution,  we observe that $V_{(g,\psi)}: \mathcal{F}_\alpha^\infty \to \mathcal{F}_\alpha^p$
is bounded if and only if
\begin{align*}
\label{carleson} \int_{\CC} |V_{(g,
\psi)}f(z)|^pe^{\frac{-p\alpha}{2}|z|^2}dm(z)\simeq \int_{\CC}
|f(z)|^p
d\mu_{(p,\alpha)}(z)\ \ \ \ \ \ \ \  \ \ \ \ \ \nonumber\\
= \int_{\CC}|f(z)|^p
e^{\frac{-p\alpha}{2}|z|^2}d\lambda_{(p,\alpha)}(z)\lesssim\|f\|_{(\infty,\alpha)}^p
\end{align*} where $d\lambda_{(p, \alpha)}(z)= e^{\frac{p\alpha}{2}|z|^2}d\mu_{(p,\alpha)}(z)$.
The inequality above  holds if and only if $\lambda_{(p,\alpha)}$ is an
$(\infty, p)$ Fock--Carleson measure. By a similar analysis, it is easy to see that $V_{(g,\psi)}: F_\alpha^\infty \to F_\alpha^p$
is compact if and only if the measure $\lambda_{(p,\alpha)}$ is an
$(\infty, p)$ vanishing Fock--Carleson measure. By Proposition~\ref{Carleson}, it follows that
the statements in (i) and (ii) are equivalent, and any one of them holds if and only
\begin{equation*}
 \widetilde{\lambda_{(P,\alpha)}}(w)= \int_{\CC} |k_{(w,\alpha)}(z)|^p
e^{\frac{-p\alpha}{2}|z|^2}d\lambda_{(p,\alpha)}(z)\in L(\CC,
dm)\end{equation*}
 from which substituting back $d\lambda_{(p,\alpha)}$ and $d\mu_{(p,\alpha)}$
 in terms of $dm,$ we obtain
\begin{equation*}
\widetilde{\lambda_{(p,\alpha)}}(w)= \int_{\CC}
|k_{(w,\alpha)}(z)|^p e^{\frac{-p\alpha}{2}|z|^2} d\lambda_{(p,
\alpha)}(z)\simeq B_{(\psi,\alpha)}(|g|^p)(w).
\end{equation*}We remain to prove the norm estimate in \eqref{normless}. But this rather can
be  deduce easily.  Since $\lambda_{(p,\alpha)}$ is an $(\infty, p)$
Fock--Carleson measure, the series of norm estimates in
Theorem~\ref{Carleson} yield
\begin{equation*}
\|V_{(g,\psi)}\|^p= \|I_{\lambda_{(p,\alpha)}}\|^p \simeq
\|\widetilde{\lambda_{(p,\alpha)}}\|_{L(\CC,dm)}\simeq
\|B_{(\psi,\alpha)}(|g|^p)\|_{L(\CC,dm)}
\end{equation*}
 and completes the proof of the theorem.\\
 \emph{Proof of Theorem~\ref{bounded4}.} Applying \eqref{comp}, we  compute
 \begin{align*}
 \|C_{(\psi,g)}f\|_{(p,\alpha)}^p \simeq \int_{\CC} |(C_{(\psi,g)}f)'(z)|^p(1+|z|)^{-p} e^{\frac{-p\alpha}{2}|z|^2}dm(z)\\
 =\int_{\CC}\frac{ |f(\psi(z)) g'(\psi(z))|^p}{(1+|z|)^p |\psi'(z)|^{-p}} e^{\frac{-p\alpha}{2}|z|^2}dm(z)\\
 =\int_{\CC}|f(\psi(z))|^p  e^{\frac{-p\alpha}{2}|\psi(z)|^2}d\upsilon(z)\\
 =\int_{\CC}|f(z)|^p  e^{\frac{-p\alpha}{2}|z|^2}d\theta(z)
  \end{align*} where $ d\theta(z)= d\upsilon(\psi^{-1}(z))$ and
  \begin{equation*}
  d \upsilon(w)= \frac{ | g'(\psi(w))|^pe^{\frac{p\alpha}{2}\big(|\psi(w)|^2-|w|^2\big)}}{(1+|w|)^p |\psi'(w)|^{-p}} dm(w).
  \end{equation*} From this it follows that the map $C_{(\psi,g)}: \mathcal{F}_\alpha^\infty \to \mathcal{F}_\alpha^p$ is bounded (compact) if and only if  $\theta$ is an $(\infty, p)$  (vanishing) Fock--Carleson measure. The remaining part of the proof  is  very similar to the  proof of the last part of  Theorem~\ref{bounded3} above and we omit it.

  We note in passing that using the above simple techniques of reformulating the bounded and compact
    properties of $C_{(\psi,g)}$ in  terms of Fock--Carleson measures,  Theorem 3.1, Theorem 3.2 and Theorem 3.3 of \cite{ZHXL}, we deduce
    the following.
  \begin{prop} Let $0<p\leq q<\infty,$ and  $\psi$ and $g$ be entire functions on  $\CC$. Then
  \begin{enumerate}
\item $C_{(\psi, g)}: \mathcal{F}_\alpha^p\to \mathcal{F}_\alpha^q$ is
     bounded if and only if
$ M_{(\psi, \alpha)}(|g|^q) \in L^\infty(\CC, dm)$, with norm estimated by
\begin{equation*}
\label{norm1estimate} \|C_{(\psi, g)}\|^q \simeq \sup_{w\in
\CC}M_{(\psi, \alpha)}(|g|^q)(w).
\end{equation*}
\item $C_{(\psi, g)}: \mathcal{F}_\alpha^p\to \mathcal{F}_\alpha^q$ compact if and only if $$\lim_{|z|\to \infty}
M_{(\psi,\alpha)}(|g|^q) (z)= 0.$$
  \item  $C_{(\psi, g)}: \mathcal{F}_\alpha^q \to \mathcal{F}_\alpha^p, \ \ \ p\neq q , $ is
  bounded (compact) if and only if $M_{(\psi,\alpha)}(|g|^p)\in L^{q/(q-p)}(\CC, dm),$ with norm estimated by ,
\begin{equation*}
\label{normless} \|C_{(\psi, g)}\|^{\frac{q}{q-p}} \simeq \int_{\CC}
M^{\frac{q}{q-p}}_{(\psi,\alpha)}(|g|^p)(w)dm(w).
\end{equation*}
    \end{enumerate}
  \end{prop}
 \section{Weighted Composition operator}
 Each pair of entire  functions $(u, \psi)$  induces a weighted composition operator defined
 by $(uC_\psi) f= u.f(\psi).$ The various mapping properties of $uC_\psi $
  acting on  a number of spaces of holomorphic  function have been studied by several authors. Its bounded and compact
properties  were
studied in \cite{SS,UEKI2} when it acts from $\mathcal{F}_\alpha^p $ to $ \mathcal{F}_\alpha^\infty $   and  $0<p\leq \infty$.  But the
problem to identify   these properties   has been open when $uC_\psi$  acts
from  $\mathcal{F}_\alpha^\infty$ to smaller spaces  $\mathcal{F}_{\alpha}^p$,\  $ 0<p<\infty$ as far as we know. In this
section,  we will answer this problem by relating it with some mapping properties of $V_{(g,\psi)}$ .

   Interestingly, the description of Fock spaces in terms of derivative, \eqref{olivia} and \eqref{compare},  makes it possible to
link some operator theoretic results on  product of  Volterra type integral and
composition operators $V_{(g,\psi)}$  with  weighted composition operators $uC_\psi$  when all acting
 between Fock spaces. For example; the bounded and compactness results on $uC_\psi: \mathcal{F}_\alpha^p \to \mathcal{F}_\alpha^\infty $
in  \cite{SS,UEKI2} could be directly read by simply replacing the weight
$|g'(z)|/(1+|z|)$  by $|u(z)|$ in the conditions of  Theorem~\ref{bounded1}. In fact, the
inequalities in \eqref{normbound}, which arise because of
\eqref{compare},  will be   replaced by equality for the case of  weighted
composition operators. By making  the same replacement in the conditions of
Theorem~\ref{bounded3}, we deduce the following result  and answer the open problem  mentioned in  the previous
paragraph.
 \begin{cor}Let $0<p< \infty, $ and  $\psi$ and $g$ be entire functions on  $\CC$.  Then the following
 statements are equivalent.
\begin{enumerate}
\item $uC_\psi: \mathcal{F}_\alpha^\infty \to \mathcal{F}_\alpha^p$ is
  bounded;
  \item $uC_\psi: \mathcal{F}_\alpha^\infty \to \mathcal{F}_\alpha^p$ is compact;
  \item  $$\int_{\CC}\int_{\CC}|k_{(w,\alpha)}(\psi(z))u(z)|^p e^{-\frac{p\alpha}{2}|z|^2} dm(z)dm(w)<\infty.$$  Moreover,  if
  $uC_\psi$ is bounded, then we have the  asymptotic estimate
  \begin{equation*}
  \label{normless}
  \|uC_\psi\|^{p} \simeq \int_{\CC}\int_{\CC}|k_{(w,\alpha)}(\psi(z))u(z)|^p e^{-\frac{p\alpha}{2}|z|^2} dm(z)dm(w)<\infty.
  \end{equation*}
  \end{enumerate}
\end{cor}
We note that  because of the derivative
on $g$ and the factor $1/(1+|z|)$, more symbols $g$ are admissible
for bounded and compact $V_{(g,\psi)}$ than  the corresponding weight $u$  for bounded
and compact $uC_\psi.$

\end{document}